\documentclass{amsart}
\usepackage{graphicx,amssymb,stmaryrd,amsfonts,epsfig,amsthm,a4,amsmath,mathrsfs,url}
\usepackage{vmargin,amscd}
\usepackage{eufrak}
\usepackage[latin1]{inputenc}
\psfigdriver{dvips}

\newtheorem{thm}{Theorem}[section]
\newtheorem{cor}[thm]{Corollary}
\newtheorem{lem}[thm]{Lemma}
\newtheorem{prop}[thm]{Proposition}
\theoremstyle{definition}
\newtheorem{defn}[thm]{Definition}

\theoremstyle{remark}
\newtheorem{obs}[thm]{Observation}
\newtheorem{rem}[thm]{Remark}

\numberwithin{equation}{section}

\newcommand{\Z}{\mathbf{Z}}

\newcommand{\R}{\mathbf{R}}
\newcommand{\C}{\mathbf{C}}
\newcommand{\Q}{\mathbf{Q}}
\newcommand{\K}{\mathbf{K}}

\newcommand{\g}{\mathfrak{g}}
\newcommand{\tpr}{\begin{tiny}\noindent Proof:}

\newcommand{\mk}{\mathfrak}
\newcommand{\Ker}{\textnormal{Ker}}

\newcommand{\SL}{\textnormal{SL}}

\newcommand{\SO}{\textnormal{SO}}
\newcommand{\SU}{\textnormal{SU}}
\newcommand{\PGL}{\textnormal{PGL}}
\newcommand{\PSL}{\textnormal{PSL}}
\newcommand{\Sp}{\textnormal{Sp}}

\newcommand{\Der}{\textnormal{Der}}

\newcommand{\Aut}{\textnormal{Aut}}
\newcommand{\Alt}{\textnormal{Alt}}
\newcommand{\Bil}{\textnormal{Bil}}

\newcommand{\rad}{\textnormal{rad}}

\newcommand{\bpr}{\noindent \textbf{Proof}: ~}

\newcommand{\epr}{~$\blacksquare$}
\begin{document}
\title{Kazhdan and Haagerup Properties in algebraic groups over local fields}
\author{Yves de Cornulier}%
\date{October 7, 2005}% (Last modification: \today)}
% --------------------------------------------------------------------------------------------------------
%                                                                                               ABSTRACT
\begin{abstract}

Given a Lie algebra $\mk{s}$, we call Lie $\mk{s}$-algebra a Lie
algebra endowed with a reductive action of $\mk{s}$. We
characterize the minimal $\mk{s}$-Lie algebras with a nontrivial
action of $\mk{s}$, in terms of irreducible representations of
$\mk{s}$ and invariant alternating forms.

As a first application, we show that if $\g$ is a Lie algebra over
a field of characteristic zero whose amenable radical is not a
direct factor, then $\g$ contains a subalgebra which is isomorphic
to the semidirect product of $\mk{sl}_2$ by either a nontrivial
irreducible representation or a Heisenberg group (this was
essentially due to Cowling, Dorofaeff, Seeger, and Wright). As a
corollary, if $G$ is an algebraic group over a local field $\K$ of
characteristic zero, and if its amenable radical is not, up to
isogeny, a direct factor, then $G(\K)$ has Property (T) relative
to a noncompact subgroup. In particular, $G(\K)$ does not have
Haagerup's property. This extends a similar result of Cherix,
Cowling and Valette for connected Lie groups, to which our method
also applies.

We give some other applications. We provide a characterization of
connected Lie groups all of whose countable subgroups have
Haagerup's property. We give an example of an arithmetic lattice
in a connected Lie group which does not have Haagerup's property,
but has no infinite subgroup with relative Property (T). We also
give a continuous family of pairwise non-isomorphic connected Lie
groups with Property (T), with pairwise non-isomorphic (resp.
isomorphic) Lie algebras.
\end{abstract}
\maketitle
% --------------------------------------------------------------------------------------------------------
%                                                                                               INTRODUCTION

\section{Introduction}

In the sequel, all Lie algebras are finite-dimensional over a
field of characteristic zero, denoted by $K$, or $\K$ when it is a
local field. If $\g$ is a Lie algebra, denote by $\rad(\g)$ its
radical and $Z(\g)$ its centre, $D\g$ its derived subalgebra, and
$\Der(\g)$ the Lie algebra of all derivations of $\g$. If
$\mk{h}_1,\mk{h}_2$ are Lie subalgebras of $\g$,
$[\mk{h}_1,\mk{h}_2]$ denotes the Lie subalgebra generated by the
brackets $[h_1,h_2]$, $(h_1,h_2)\in\mk{h}_1\times\mk{h}_2$.

Let $\g$ be a Lie algebra with radical $\rad(\g)=\mk{r}$ and
semisimple Levi factor $\mk{s}$ (so that $\g\simeq\mk{s}\ltimes
\mk{r}$). We focus here on aspects of $\g$ related to the action
of $\mk{s}$. This suggests the following definitions.

If $\mk{s}$ is a Lie algebra, we define a {\em Lie
$\mk{s}$-algebra} to be a Lie algebra $\mk{n}$ endowed with a
morphism $i:\mk{s}\to \Der(\mk{n})$, defining a {\em completely
reducible} action of $\mk{s}$ on $\mk{n}$. (This latter technical
condition is empty if $\mk{s}$ is semisimple.)

A Lie $\mk{s}$-algebra naturally embeds in the semidirect product
$\mk{s}\ltimes\mk{n}$, so that we write $i(s)(n)=[s,n]$ for
$s\in\mk{s}$, $n\in\mk{n}$.

By the {\em trivial irreducible module} of $\mk{s}$ we mean a
one-dimensional vector space endowed with a trivial action of
$\mk{s}$. We say that a module (over a Lie algebra or over a
group) is {\em full} if is is completely reducible and does not
contain the trivial irreducible module.

\begin{defn} Let $\mk{s}$ be a Lie algebra. We say that
a Lie $\mk{s}$-algebra $\mk{n}$ is {\em minimal} if
$[\mk{s},\mk{n}]\neq 0$, and for every $\mk{s}$-subalgebra
$\mk{n}'$ of $\mk{n}$, either $\mk{n}'=\mk{n}$ or
$[\mk{s},\mk{n}']=0$.
\end{defn}

It is clear that a Lie $\mk{s}$-algebra $\mk{n}$ satisfying
$[\mk{s},\mk{n}]\neq 0$ contains a minimal $\mk{s}$-subalgebra. We
begin by a characterization of minimal $\mk{s}$-algebras:

\begin{thm}
Let $\mk{s}$ be a Lie algebra. A solvable Lie $\mk{s}$-algebra
$\mk{n}$ is minimal if and only if it satisfies the following
conditions 1), 2), 3), and 4):

1) $\mk{n}$ is 2-nilpotent (that is, $[\mk{n},D\mk{n}]=0$).

2) $[\mk{s},\mk{n}]=\mk{n}$.

3) $[\mk{s},D\mk{n}]=0$.

4) $\mk{n}/D\mk{n}$ is irreducible as a
$\mk{s}$-module.\label{car_minimal}
\end{thm}

\begin{defn}
We call a solvable Lie $\mk{s}$-algebra $\mk{n}$ {\em almost
minimal} if it satisfies conditions 1), 2), and 3) of Theorem
\ref{car_minimal}.\label{def_alm_min}
\end{defn}
This definition has the advantage to be invariant under field
extensions. Note that an almost minimal solvable Lie
$\mk{s}$-algebra $\mk{n}$ automatically satisfies the following
Condition 4'): $\mk{n}/D\mk{n}$ is a full $\mk{s}$-module.

\medskip

The classification of (almost) minimal solvable Lie
$\mk{s}$-algebras can be deduced from the classification of
irreducible $\mk{s}$-modules. Let $\mk{v}$ be a full
$\mk{s}$-module (equivalently, an abelian Lie $\mk{s}$-algebra
satisfying $[\mk{s},\mk{v}]=\mk{v}$). Recall that a bilinear form
$\varphi$ on $\mk{v}$ is called $\mk{s}$-invariant if it satisfies
$\varphi([s,v],w)+\varphi(v,[s,w])=0$ for all $s\in\mk{s}$,
$v,w\in\mk{v}$. Let $\Bil_\mk{s}(\mk{v})$ (resp.
$\Alt_{\mk{s}}(\mk{v})$) denote the space of all
$\mk{s}$-invariant bilinear (resp. alternating bilinear) forms on
$\mk{v}$. Denote by $\Alt_{\mk{s}}(\mk{v})^*$ the linear dual
of~$\Alt_{\mk{s}}(\mk{v})$.

\begin{defn}
We define the Lie $\mk{s}$-algebra $\mk{h}(\mk{v})$ as follows: as
a vector space, $\mk{h}(\mk{v})=\mk{v}\oplus
\Alt_{\mk{s}}(\mk{v})^*$; it is endowed with the following
bracket:
\begin{equation}[(x,z),(x',z')]=(0,e_{x,x'})\qquad x,x'\in\mk{v}\;\;
z,z'\in \Alt_{\mk{s}}(\mk{v})^*\label{LP}\end{equation} where
$e_{x,x'}\in \Alt_{\mk{s}}(\mk{v})^*$ is defined by
$e_{x,x'}(\varphi)=\varphi(x,x')$.\label{hvZ}\end{defn}

This is a 2-nilpotent Lie $\mk{s}$-algebra under the action
$[s,(x,z)]=([s,x],0)$, which is almost minimal. Other almost
minimal Lie $\mk{s}$-algebras can be obtained by taking the
quotient by a linear subspace of the centre. The following theorem
states that this is the only way to construct almost minimal
solvable Lie $\mk{s}$-algebras.

\begin{thm}
If $\mk{n}$ is an almost minimal solvable Lie $\mk{s}$-algebra,
then it is isomorphic (as a $\mk{s}$-algebra) to
$\mk{h}(\mk{v})/Z$, for some full $\mk{s}$-module $\mk{v}$ and
some subspace $Z$ of $\Alt_{\mk{s}}(\mk{v})^*$. It is minimal if
and only if $\mk{v}$ is irreducible.

Moreover, the almost minimal $\mk{s}$-algebras $\mk{h}(\mk{v})/Z$
and $\mk{h}(\mk{v})/Z'$ are isomorphic if and only if $Z'$ and $Z$
are in the same orbit for the natural action of
$\Aut_{\mk{s}}(\mk{v})$ on the Grassmannian of
$\Alt_{\mk{s}}(\mk{v})^*$.\label{thm:clasmin}
\end{thm}

\begin{rem}
If $\mk{s}$ is semisimple, $\mk{s}\ltimes \mk{h}(\mk{v})$ is the
universal central extension of the perfect Lie algebra
$\mk{s}\ltimes \mk{v}$.
\end{rem}

The case of $\mk{sl}_2$ is essential, and there is a simple
description for it. Recall that if $\mk{s}=\mk{sl}_2$, then, up to
isomorphism, there exists exactly one irreducible $\mk{s}$-module
$\mk{v}_{n}$ of dimension $n$ for every $n\ge 1$. If $n=2m$ is
even, it has a central extension by a one-dimensional subspace,
giving a Heisenberg Lie algebra, on which $\mk{sl}_2$ acts
naturally (see \ref{section_sl2} for details), denoted by
$\mk{h}_{2m+1}$. Theorem \ref{thm:clasmin} thus reduces as:

\begin{prop}
Up to isomorphism, the minimal solvable Lie $\mk{sl}_2$-algebras
are $\mk{v}_n$ and $\mk{h}_{2n-1}$ ($n\ge 2$).\label{sl2min}
\end{prop}

Let $\g$ be a Lie algebra, $\mk{r}$ its radical and $\mk{s}$ a
semisimple factor. Write $\mk{s}=\mk{s}_{c}\oplus\mk{s}_{nc}$ by
separating anisotropic and isotropic factors\footnote{$c$ and $nc$
respectively stand for ``non-compact" and ``compact"; this is
related to the fact that if $S$ is a simple algebraic group
defined over the local field $\K$, then its Lie algebra is
$\K$-isotropic if and only if $S(\K)$ is not compact.}. The ideal
$\mk{s}_c\ltimes\mk{r}$ is sometimes called the {\em amenable
radical} of $\g$.

\begin{defn}
We call $\g$ M-decomposed if $[\mk{s}_{nc},\mk{r}]=0$.
Equivalently, $\g$ is M-decomposed if the amenable radical is a
direct factor of $\g$.
\end{defn}

\begin{prop}
Let $\g$ be a Lie algebra, and keep notation as above. Suppose
that $\g$ is not M-decomposed. Then there exists a Lie subalgebra
$\mk{h}$ of $\g$ which is isomorphic to $\mk{sl}_2\ltimes\mk{v}_n$
or $\mk{sl}_2\ltimes\mk{h}_{2n-1}$ for some $n\ge
2$.\label{subalg_snc}
\end{prop}

This result is essentially due to \cite{CDSW}, where it is not
explicitly stated, but it is actually proved in the proof of
Proposition 8.2 there (under the assumption $K=\R$, but their
argument generalizes to any field of characteristic zero). This
was a starting point for the present paper.

\medskip

Let $G$ be a locally compact, $\sigma$-compact group. Recall that
$G$ has the Haagerup Property if it has a metrically proper
isometric action on a Hilbert space; in contrast, $G$ has
Kazhdan's Property (T) if every isometric action of $G$ on a
Hilbert space has a fixed point. See \ref{reminder} for a short
reminder about Haagerup and Kazhdan Properties.

\medskip

We provide corresponding statements for Proposition
\ref{subalg_snc} in the realm of algebraic groups and connected
Lie groups. As a consequence, we get the following theorem, which
was the initial motivation for the results above. It was already
proved, in a different way, for connected Lie groups in
\cite[Chap. 4]{CCJJV}.

\begin{thm}
Let $G$ be either a connected Lie group, or $G=\mathbf{G}(\K)$,
where $\mathbf{G}$ is a linear algebraic group over the local
field $\K$ of characteristic zero. Let $\g$ be its Lie algebra.
The following are equivalent.

\begin{itemize}
\item[(i)] $G$ has Haagerup's property.

\item[(ii)] For every noncompact closed subgroup $H$ of $G$,
$(G,H)$ does not have relative Property~(T).

\item[(iii)] The following conditions are satisfied:
\begin{itemize}
\item $\mathfrak{g}$ is M-decomposed. \item All simple factors of
$\mathfrak{g}$ have $\mathbf{K}$-rank~$\le 1$. \item (in the case
of Lie groups or when $\mathbf{K}=\mathbf{R}$) No simple factor of
$\mathfrak{g}$ is isomorphic to $\mathfrak{sp}(n,1)$ ($n\ge 2$)
or~$\mathfrak{f}_{4(-20)}$.
\end{itemize}
\item[(iv)] $\mathfrak{g}$ contains no isomorphic copy of any one
of the following Lie algebras
\begin{itemize}
\item $\mathfrak{sl}_2\ltimes\mathfrak{v}_n$ or
$\mathfrak{sl}_2\ltimes\mathfrak{h}_{2n-1}$ for some $n\ge 2$,
\item (in the case of Lie groups or when $\mathbf{K}=\mathbf{R}$)
$\mathfrak{sp}(2,1)$.\end{itemize}
\end{itemize}\label{class_haag}
\end{thm}

\begin{rem}
The notion of M-decomposed (real) Lie algebras also appears in
other contexts: heat kernel on Lie groups \cite{Var}, Rapid Decay
Property \cite{CPS}, weak amenability \cite{CDSW}.
\end{rem}

\bigskip

We derive some other results with the help of Theorem
\ref{thm:clasmin}.

\begin{prop}
There exists a continuous family $(\g_t)$ of pairwise
non-isomorphic real (or complex) Lie algebras satisfying the
following properties:

(i) $\g_t$ is perfect, and

(ii) the simply connected Lie group corresponding to $\g_t$ has
Property (T).\label{cont_family_T1}
\end{prop}

Note that Proposition \ref{cont_family_T1} with only (i) may be of
independent interest; we do not know if it had already been
observed. On the other hand, it is well-known that there exist
continuously many pairwise non-isomorphic complex n-dimensional
nilpotent Lie algebras if $n\ge 7$.

\begin{prop}
There exists a continuous family of pairwise non-isomorphic
connected Lie groups with Property (T), and with isomorphic Lie
algebras.\label{cont_family_T2}
\end{prop}

We also give the classification, when $\K=\R$, of the minimal
$\mk{so}_3$-algebras (Proposition \ref{so3_min}). We use it to
prove (ii)$\Rightarrow$(i) in the following result (while the
reverse implication is essentially due to \cite[Theorem
5.1]{GHW}).

\begin{thm}
Let $G$ be a connected Lie group. Then the following are equivalent:\\
(i) $G$ is locally isomorphic to $\SO_3(\R)^a\times
\SL_2(\R)^b\times
\SL_2(\C)^c\times R$, for a solvable Lie group $R$ and integers $a,b,c$.\\
(ii) Every countable subgroup of $G$ has Haagerup's property (when
endowed with the discrete topology).\label{GdHaagerup}
\end{thm}

\begin{rem}
Assertion (i) of Theorem \ref{GdHaagerup} is equivalent to: (ii')
{\em The complexification $\g_\C$ of $\g$ is M-decomposed, and its
semisimple part is isomorphic to $\mk{sl}_2(\C)^n$ for some $n$}.
\end{rem}

For instance, $\SO_3(\R)\ltimes\R^3$ has a countable subgroup
which does not have Haagerup's property. An explicit example is
given by $\SO_3(\Z[1/p])\ltimes\Z[1/p]^3$. It can also be shown
that this group has no infinite subgroup with relative Property
(T). This answers an open question in \cite[Section 7.1]{CCJJV}.
This group is not finitely presented (this is a consequence of
\cite[Theorem 2.6.4]{Abels}); we give a similar example in Remark
\ref{SansTrel_niHaag} which is, in addition, finitely presented.

%%%%%%%%%%%%%%%%%%%%%%%%%%%%%%%%%%%%%%%%%%%%%%%%%%%%%%%%%%%%%%%%%%%%%%%%%%%%%%%%%%%%%%%%%%%%%%%%%%%%%%

% --------------------------------------------------------------------------------------------------------
%                                                                           section{minimal subalgebras}
%                                                                  subsection{proof main results}

\section{Lie algebras}

\subsection{Minimal subalgebras}

\begin{prop}
Let $\mk{n}$ be a solvable Lie $\mk{s}$-algebra.

1) The Lie $\mk{s}$-subalgebra $[\mk{s},\mk{n}]$ is an ideal in
$\mk{n}$ (and also in $\mk{s}\ltimes\mk{n}$), and
$[\mk{s},[\mk{s},\mk{n}]]=[\mk{s},\mk{n}]$.

2) If, moreover, $[\mk{s},D\mk{n}]=0$, then $[\mk{s},\mk{n}]$ is
an almost minimal Lie algebra (see Definition
\ref{def_alm_min}).\label{[s,n]almost_minimal}
\end{prop}
\bpr 1) Let $\mk{v}$ be the subspace generated by the brackets
$[s,n]$, $(s,n)\in\mk{s}\times\mk{n}$. Since the action of
$\mk{s}$ is completely reducible (see the definition of Lie
$\mk{s}$-algebra), it is immediate that $[\mk{s},\mk{n}]$ and
$[\mk{s},[\mk{s},\mk{n}]]$ both coincide with the Lie subalgebra
generated by $\mk{v}$. Then, using Jacobi identity,
$$[\mk{n},[\mk{s},\mk{n}]]\;=\;[\mk{n},[\mk{s},[\mk{s},\mk{n}]]]\;\subseteq\;
[\mk{s},[\mk{n},[\mk{s},\mk{n}]]]\;+\;[[\mk{s},\mk{n}],[\mk{s},\mk{n}]]\;\subseteq\;
[\mk{s},[\mk{n},\mk{n}]]\;+\;[\mk{s},\mk{n}]\;\subseteq\;
[\mk{s},\mk{n}].$$

\medskip

2) Let $\mk{z}$ be the linear subspace generated by the
commutators $[v,w]$, $v,w\in\mk{v}$. By Jacobi identity,
$$[\mk{v},\mk{z}]\;=\;[[\mk{s},\mk{v}],\mk{z}]\;\subseteq\;
[[\mk{s},\mk{z}],\mk{v}]+[\mk{s},[\mk{v},\mk{z}]]\;\subseteq\;
[[\mk{s},D\mk{n}],\mk{v}]+[\mk{s},D\mk{n}]=0.$$

Thus, the subspace $\mk{n}'=\mk{v}\oplus \mk{z}$ is a 2-nilpotent
Lie $\mk{s}$-subalgebra of $\mk{n}$. The Lie subalgebra
$[\mk{s},\mk{n}']$ contains $\mk{v}$, hence also contains
$\mk{z}$, so $[\mk{s},\mk{n}']$ is equal to $\mk{n}'$. Thus
Conditions 1) and 2) of Definition \ref{def_alm_min} are
satisfied, while Condition 3) follows immediately from the
hypothesis $[\mk{s},D\mk{n}]=0$.\epr

\medskip

%\begin{thm}
%Let $\mk{s}$ be a Lie algebra, and $\mk{n}$ a solvable Lie
%$\mk{s}$-algebra. The following are equivalent.

%(i) $\mk{n}$ is a minimal Lie $\mk{s}$-algebra.

%(ii) The Lie $\mk{s}$-algebra $\mk{n}$ satisfies the following
%conditions 1), 2), 3), and 4):
%\begin{itemize}
%\item 1) $\mk{n}$ is 2-nilpotent (that is, $[\mk{n},D\mk{n}]=0$).

%\item 2) $[\mk{s},\mk{n}]=\mk{n}$.

%\item 3) $[\mk{s},D\mk{n}]=0$.

%\item 4) $\mk{n}/D\mk{n}$ is irreducible as a $\mk{s}$-module.
%\end{itemize}\label{car_minimal}
%\end{thm}
%\bpr

%\newcommand{\bprr}{\noindent \textbf{Proof} }

\noindent \textbf{Proof of Theorem \ref{car_minimal}.} Suppose
that the four conditions of the theorem are satisfied. Condition 4
implies $\mk{n}\neq 0$. Then Condition 2 implies
$[\mk{s},\mk{n}]=\mk{n}\neq 0$. Let $\mk{n}'\subseteq\mk{n}$ be a
$\mk{s}$-subalgebra. Then, by irreducibility (Condition 4), either
$D\mk{n}+\mk{n}'=D\mk{n}$ or $D\mk{n}+\mk{n}'=\mk{n}$. In the
first case, $\mk{n}'$ centralizes $\mk{s}$. In the second case,
$\mk{n}=[\mk{s},\mk{n}]=[\mk{s},\mk{n'}+D\mk{n}]=[\mk{s},\mk{n}']\subseteq\mk{n'}$,
using Conditions 1 and 2, and the fact that $\mk{n}'$ is a
$\mk{s}$-subalgebra.

\medskip

Conversely, suppose that the $\mathfrak{s}$-algebra $\mathfrak{n}$
is solvable and minimal. Since $\mk{n}$ is solvable, $D\mk{n}$ is
a proper $\mk{s}$-subalgebra, so that, by minimality,
$[\mk{s},D\mk{n}]=0$. By Proposition \ref{[s,n]almost_minimal},
$[\mk{s},\mk{n}]$ is a nonzero almost minimal Lie
$\mk{s}$-subalgebra of $\mk{n}$, hence satisfies 1), 2), 3). The
minimality implies that 4) is also satisfied.\epr

%\begin{thm}
%If $\mk{n}$ is an almost minimal solvable Lie $\mk{s}$-algebra,
%then it is isomorphic (as a $\mk{s}$-algebra) to
%$\mk{h}(\mk{v})/Z$, for some full $\mk{s}$-module $\mk{v}$ and
%some subspace $Z$ of $\Alt_{\mk{s}}(\mk{v})^*$. It is minimal if
%and only if $\mk{v}$ is irreducible.

%Moreover, the almost minimal $\mk{s}$-algebras $\mk{h}(\mk{v})/Z$
%and $\mk{h}(\mk{v})/Z'$ are isomorphic if and only if $Z'$ and $Z$
%are in the same orbit for the natural action of
%$\Aut_{\mk{s}}(\mk{v})$ on the Grassmannian of
%$\Alt_{\mk{s}}(\mk{v})^*$.label{thm:clasmin}
%\end{thm}

\medskip

\noindent \textbf{Proof of Theorem \ref{thm:clasmin}.} Let
$\mk{n}$ be an almost minimal solvable Lie $\mk{s}$-algebra. Let
$\mk{v}$ be the subspace generated by the brackets $[s,n]$,
$(s,n)\in\mk{s}\times\mk{n}$. Since $\mk{n}$ is almost minimal,
$\mk{v}$ is a complementary subspace of $D\mk{n}$, and is a full
$\mk{s}$-module. If $u\in D\mk{n}^*$, consider the alternating
bilinear form $\phi_u$ on $\mk{v}$ defined by
$\phi_u(x,y)=u([x,y])$. This defines a mapping
$D\mk{n}^*\to\Alt_{\mk{s}}(\mk{v})$ which is immediately seen to
be injective. By duality, this defines a surjective linear map
$\Alt_{\mk{s}}(\mk{v})^*\to D\mk{n}$, whose kernel we denote by
$Z$. It is immediate from the definition of $\mk{h}(\mk{v})$ that
this map extends to a surjective morphism of Lie $\mk{s}$-algebras
$\mk{h}(\mk{v})\to\mk{n}$ with kernel $Z$. This proves that
$\mk{n}$ is isomorphic to $\mk{h}(\mk{v})/Z$.

The second assertion is immediate.

The third assertion follows from the proof of the first one, where
we made no choice. Namely, take an isomorphism
$\psi:\mk{h}(\mk{v})/Z\to\mk{h}(\mk{v})/Z'$. It gives by
restriction an $\mk{s}$-automorphism $\varphi$ of $\mk{v}$, which
induces a unique automorphism $\tilde{\varphi}$ of
$\mk{h}(\mk{v})$. Let $p$ and $p'$ denote the natural projections
in the following diagram of Lie $\mk{s}$-algebras:
\[
\begin{CD}
{\mk{h}(\mk{v})} @> {p} >> {\mk{h}(\mk{v})/Z} \\
@ V{\tilde{\varphi}} VV @ VV {\psi} V \\
{\mk{h}(\mk{v})} @>> {p'} > {\mk{h}(\mk{v})/Z'} \\
\end{CD}
\]
This diagram is commutative: indeed, $p'\circ\tilde{\varphi}$ and
$\psi\circ p$ coincide in restriction to $\mk{v}$, and $\mk{v}$
generates $\mk{h}(\mk{v})$ as a Lie algebra. This implies
$Z=\Ker(\psi\circ
p)=\Ker(p'\circ\tilde{\varphi})=\tilde{\varphi}^{-1}(Z')$.\epr

\medskip

% --------------------------------------------------------------------------------------------------------
%                                                                           subsection{example sl2}

\subsection{The example $\mk{sl}_2$}\label{section_sl2}

If $\mk{s}=\mk{sl}_2(K)$, then, up to isomorphism, there exists
exactly one irreducible $\mk{s}$-module $\mk{v}_{n}$ of dimension
$n$ for every $n\ge 1$.

Since $\mk{v}_{n}$ is absolutely irreducible for all $n$, by
Schur's Lemma, $\Bil_\mk{s}(\mk{v}_{n})$ is at most one
dimensional for all $n$. In fact, it is one-dimensional. Indeed,
take the usual basis $(H,X,Y)$ of $\mk{sl}_2$ satisfying
$[H,X]=2X$, $[H,Y]=-2Y$, $[X,Y]=H$, and take the basis
$(e_0,\dots,e_{n-1})$ of $\mk{v}_{n}$ so that $H.e_i=(n-1-2i)e_i$,
$X.e_i=(n-i)e_{i-1}$, and $Y.e_i=(i+1)e_{i+1}$, with the
convention $e_{-1}=e_{n}=0$. Then $Bil_\mk{s}(\mk{v}_{n})$ is
generated by the form $\varphi_{n}$ defined by
$$\varphi_n(e_i,e_{n-1-i})=(-1)^i\begin{pmatrix}
  i \\
  n-1 \\
\end{pmatrix};\quad\varphi(e_i,e_j)=0 \textnormal{ if } i+j\neq
n-1.$$ For odd $n$, $\varphi_{n}$ is symmetric so that
$Alt_\mk{s}(\mk{v}_{n})=0$; for even $n$, $\varphi_{n}$ is
symplectic and generates $\Alt_\mk{s}(\mk{v}_{n})$. For even $n$,
denote by $\mk{h}_{n+1}$ the one-dimensional central extension
$\mk{h}(\mk{v}_n)$, well-known as the $(n+1)$-dimensional
Heisenberg Lie algebra.

\medskip

\noindent \textbf{Proof of Proposition \ref{subalg_snc}.} Since
$\mk{s}_{nc}$ is semisimple and isotropic, it is generated by its
subalgebras $K$-isomorphic to $\mk{sl}_2$. Since
$[\mk{s}_{nc},\mk{r}]\neq 0$, this implies that there exists some
subalgebra $\mk{s}'$ of $\mk{s}_{nc}$ which is $K$-isomorphic to
$\mk{sl}_2$ and such that $[\mk{s}',\mk{r}]\neq 0$. Then the
result is clear from Proposition \ref{sl2min}. Notice that the
proof gives the following slight refinement: $\mk{h}$ can be
chosen so that $\rad(\mk{h})\subseteq\rad(\g)$.\epr

\medskip

% --------------------------------------------------------------------------------------------------------
%                                                                           section{ex so3}

\subsection{The example $\mk{so}_3$}\label{section_so3}

We now study a more specific example. Let us deal with the field
$\R$ of real numbers, and with $\mk{s}=\mk{so}_3$.

Since the complexification of $\mk{so}_3$ is isomorphic to
$\mk{sl}_2(\C)$, the irreducible complex $\mk{s}$-modules make up
a family $(\mk{d}_n^\C)$ $(n\ge 1)$; $\dim_\C(\mk{d}_n^\C)=n$,
which are the symmetric powers of the natural action of
$\mk{su}_2=\mk{so}_3 $ on $\C^2$.

If $n=2m+1$ is odd, then this is the complexification of a real
$\mk{so}_3$-module $\mk{d}_{2m+1}^\R$ (of dimension $n$). If
$n=2m$ is even, $\mk{d}_n^\C$ is irreducible as a $4m$-dimensional
real $\mk{so}_3$-module, we call it $\mk{u}_{4m}$.

These two families $(\mk{d}_{2n+1}^\R)$ and $(\mk{u}_{4n})$ make
up all irreducible real $\mk{so}_3$-modules.

\begin{prop}
The irreducible real $\mk{so}_3$-modules make up two families: a
family $(\mk{d}_{2n+1}^\R)$ of $(2n+1)$-dimensional modules $(n\ge
0$), absolutely irreducible, and a family $(\mk{u}_{4n})$ of
$4n$-dimensional modules ($n\ge 1$), not absolutely irreducible,
preserving a quaternionic structure.\epr\label{repsu2}
\end{prop}

Since $(\mk{d}_{2n+1}^\R)$ is absolutely irreducible, the space of
invariant bilinear forms on $(\mk{d}_{2n+1}^\R)$ is generated by a
scalar product, so that $\Alt_{\mk{so}_3}(\mk{d}_{2n+1}^\R)=0$

On the other hand, $\Alt_{\mk{so}_3}(\mk{u}_{4n})$ is
three-dimensional, and is given by the imaginary part of an
invariant quaternionic hermitian form.

In order to classify the minimal solvable $\mk{so}_3$-algebras, we
must determine the orbits of the natural action of
$\Aut_{\mk{so}_3}(\mk{u}_{4n})$ on
$\Alt_{\mk{so}_3}(\mk{u}_{4n})$. It is a standard fact that
$\Aut_{\mk{so}_3}(\mk{u}_{4n})$ is isomorphic to the group of
nonzero quaternions, that $\Alt_{\mk{so}_3}(\mk{u}_{4n})$
naturally identifies with the set of imaginary quaternions, and
that the action of $\Aut_{\mk{so}_3}(\mk{u}_{4n})$ on
$\Alt_{\mk{so}_3}(\mk{u}_{4n})$ is given by conjugation of
quaternions. This implies that it acts transitively on each
component of the Grassmannian.

For $i=0,1,2,3$, let $Z_i$ be a fixed $(3-i)$-dimensional linear
subspace of $\Alt_{\mk{s}}(\mk{v})^*$. Denote by $\mk{hu}_{4n}^i$
the minimal Lie $\mk{so}_3$-algebra $\mk{h}(\mk{u}_{4n})/Z_i$; of
course, $\mk{hu}_{4n}^0=\mk{u}_{4n}$ and
$\mk{hu}_{4n}^3=\mk{h}(\mk{u}_{4n})$.

\begin{prop}
Up to isomorphism, the minimal solvable Lie
$\mk{so}_3(\R)$-algebras are $\mk{d}_{2n+1}^\R$ ($n\ge 1$) and
$\mk{hu}_{4n}^i$ ($n\ge 1$, $i=0,1,2,3$).\epr\label{so3_min}
\end{prop}

There is an analogous statement to Proposition \ref{subalg_snc}.

\begin{prop}
Let $\g$ be a Lie algebra over $\R$. Suppose that
$[\mk{s}_{c},\mk{r}]\neq\{0\}$. Then $\g$ has a Lie subalgebra
which is isomorphic to either $\mk{so}_3\ltimes \mk{d}_{2n+1}^\R$
or $\mk{so}_3\ltimes \mk{hu}_{4n}^i$ for some $i=0,1,2,3$ and some
$n\ge 1$.\epr\label{so3sg}
\end{prop}

% --------------------------------------------------------------------------------------------------------
%                                                                section{corresponding results...}
%                                                                      subsection{minimal algebraic S-groups}

\section{Corresponding results for algebraic groups and connected Lie
groups}

\subsection{Minimal algebraic subgroups}

We now give the corresponding statements and results for algebraic
groups.

Let $S$ be a reductive $K$-group. A $K$-$S$-group means a linear
$K$-group endowed with a $K$-action of $S$ by automorphisms.

Recall that the Lie algebra functor gives an equivalence of
categories between the category of unipotent $K$-groups and the
category of nilpotent Lie $K$-algebras. If $S$ is semisimple and
simply connected with Lie algebra $\mk{s}$, it induces an
equivalence of categories between the category of unipotent
$K$-$S$-groups and the category of nilpotent Lie $S$-algebras over
$K$. If $S$ is not simply connected (in particular, if $S$ is not
semisimple), this is no longer an essentially surjective functor,
but it remains fully faithful.

A minimal (resp. almost minimal) solvable $S$-group $N$ is defined
similarly as in the case of Lie algebras; it is automatically
unipotent (since it satisfies $[S,N]=N$). Moreover, $N$ is a
minimal (resp. almost minimal) solvable $K$-$S$-group if and only
if its Lie algebra $\mk{n}$ is a minimal (resp. almost minimal)
solvable Lie $\mk{s}$-algebra. Proposition
\ref{[s,n]almost_minimal} and Theorem \ref{car_minimal} also
immediately carry over into the context of algebraic groups.
\medskip

If $S$ is reductive and $V$ is a $K$-$S$-module, we define the
unipotent $K$-$S$-group $H(V)$ as follows: as a variety,
$H(V)=V\oplus \Alt_{S}(V)^*$; it is endowed with the following
group law:
\begin{equation}(x,z)(x',z')=(x+x',z+z'+e_{x,x'})\qquad x,x'\in V\;\;
z,z'\in \Alt_{S}(V)^*\label{LPgr}\end{equation} where $e_{x,x'}\in
\Alt_{S}(V)^*$ is defined by $e_{x,x'}(\varphi)=\varphi(x,x')$.
This is a $K$-$S$-group under the action $s.(x,z)=(s.x,z)$. It is
clear that its Lie algebra is isomorphic as a Lie $K$-$S$-algebra
to $\mk{h}(\mk{v})$, where $V=\mk{v}$ viewed as a $\mk{s}$-module.
Here is the analog of Theorem \ref{thm:clasmin}.

\begin{thm}
If $N$ is an almost minimal solvable $K$-$S$-group, then it is
isomorphic (as a $K$-$S$-group) to $H(V)/Z$, for some full
$K$-$S$-module $V$ and some $K$-subspace $Z$ of $\Alt_{S}(V)^*$.
It is minimal if and only if $V$ is irreducible.

Moreover, the almost minimal $K$-$S$-groups $H(V)/Z$ and $H(V)/Z'$
are isomorphic if and only if $Z'$ and $Z$ are in the same orbit
for the natural action of $\Aut_{S}(V)$ on the Grassmannian of
$\Alt_{S}(V)^*$.\label{clasmin}
\end{thm}

% -----------------------------------------------------------------

\medskip

% --------------------------------------------------------------------------------------------------------
%                                                                           subsection{example SL2}

\subsection{The example $\SL_2$}\label{section_SL2}

The simply connected $K$-group with Lie algebra $\mk{sl}_2$ is
$\SL_2$. Denote by $V_n$ and $H_{2n-1}$ the $\SL_2$-groups
corresponding to $\mk{v}_n$ and $\mk{h}_{2n-1}$. These are the
solvable minimal $\SL_2$-groups over $K$. The only non-simply
connected $K$-group with Lie algebra $\mk{sl}_2$ is the adjoint
group $\PGL_2$; thus the minimal solvable $\PGL_2$-groups over $K$
are $V_{2n-1}$ for $n\ge 2$.

\begin{rem}It is convenient, in algebraic groups, to deal with the
unipotent radical rather than with the radical. It is
straightforward to see that a reductive subgroup $S$ of a linear
algebraic group centralizes the radical if and only if it
centralizes the unipotent radical. Indeed, suppose $[S,R_u]=1$. We
always have $[S,R/R_u]=1$ since $R/R_u$ is central in $G^0/R_u$
and $S$ is connected ($G^0$ denoting the unit component of $G$).
Since $S$ is reductive, this implies that $S$ acts
trivially on~$R$.%\footnote{Write, for $s\in S$ and $r\in R$,
%$s.r=ru(s,r)$, where $u(s,r)\in R_u$ and $u(s,r)=1$ if $r\in R_u$.
%Then, for all $s,t\in S$ and $r\in R$
%$st.r=s.ru(t,r)=(s.r)(s.u(t,r))=ru(s,r)u(t,r)$, so that
%$u(st,r)=u(s,r)u(t,r)$. This implies that if $s\in D^n S=S$, then
%$u(s,r)\in D^nR_u$. Taking $n$ sufficiently large, we obtain
%$u(s,r)=1$ for all $s\in S$ and $r\in R$, that is, $[S,R]=1$.}.
\end{rem}

Let $G$ be a linear algebraic group over the field $K$ of
characteristic zero, $R$ its radical, $S$ a Levi factor,
decomposed as $S_{nc}S_{c}$ by separating $K$-isotropic and
$K$-anisotropic factors.

\begin{prop}
Suppose that $[S_{nc},R]\neq 1$. Then $G$ has a $K$-subgroup which
is $K$-isomorphic to either $\SL_2\ltimes V_n$, $\PGL_2\ltimes
V_{2n-1}$, or $\SL_2\ltimes H_{2n-1}$ for some $n\ge
2$.\label{SL2min}
\end{prop}

Let us mention the translation into the context of connected Lie
groups, which is immediate from the Lie algebraic version.

\begin{prop}
Let $G$ be a real Lie group. Suppose that $[S_{nc},R]\neq 1$. Then
there exists a Lie subgroup $H$ of $G$ which is locally isomorphic
to $\SL_2(\R)\ltimes V_n(\R)$ or $\SL_2(\R)\ltimes H_{2n-1}(\R)$
for some $n\ge 2$.\label{subliegrp_snc}
\end{prop}

\begin{rem}
1) An analogous result holds with complex Lie groups.

2) The Lie subgroup $H$ is not necessarily closed; this is due to
the fact that $\widetilde{\SL_2(\R)}$ and $H_{2n-1}(\R)$ have
noncompact centre. For instance, take an element $z$ of the centre
of $H$ that generates an infinite discrete subgroup, and take the
image of $H$ in the quotient of $H\times\R/\Z$ by $(z,\alpha)$,
where $\alpha$ is irrational.

3) It can be easily be shown that, if the Lie group $G$ is linear,
then the subgroup $H$ is necessarily closed. In a few words, this
is because the derived subgroup of the radical is unipotent, hence
simply connected, and the centre of the semisimple part is finite.
\end{rem}

\medskip

% --------------------------------------------------------------------------------------------------------
%                                                                           section{ex SO3}

\subsection{The example $\SO_3$} We go on with the notation of
\ref{section_so3}.

In the context of algebraic $\R$-groups as in the context of
connected Lie groups, the simply connected group corresponding to
$\mk{so}_3(\R)$ is $\SU(2)$. The only non-simply connected
corresponding group is $\SO_3(\R)$.

The irreducible $\SU(2)$-modules corresponding to
$\mk{d}_{2m+1}^\R$ and $\mk{u}_{4n}$ are denoted by $D_{2n+1}^\R$
and $U_{4n}$. Among those, only $D_{2n+1}^\R$ provide
$\SO_3(\R)$-modules.

Denote by $HU_{4n}^i$ the unipotent $\R$-group corresponding to
$\mk{hu}_{4n}^i$, $i=0,1,2,3$.

\begin{rem}
It can be shown that the maximal unipotent subgroups of $\Sp(n,1)$
are isomorphic to $HU_{4n}^3$.
\end{rem}

\begin{prop}
Up to isomorphism, the minimal solvable Lie $\SO_3(\R)$-algebras
are $D_{2n+1}^\R$ for $n\ge 1$; the other minimal solvable Lie
$\SU(2)$-algebras are $HU_{4n}^i$, for $n\ge 1$, $i=0,1,2,3$.\epr
\end{prop}

\begin{prop}
Let $G$ be a linear algebraic $\R$-group. Suppose that
$[S_{c},R]\neq 1$. Then $G$ has a $\R$-subgroup which is
$\R$-isomorphic to either $\SU(2)\ltimes D_{2n+1}^\R$,
$\SO_3(\R)\ltimes D_{2n+1}^\R$, or $\SU(2)\ltimes HU_{4n}^i$ for
some $i=0,1,2,3$ and some $n\ge 1$.

Let $G$ be a real Lie group. Suppose that $[S_{c},R]\neq 1$. Then
$G$ has a Lie subgroup which is locally isomorphic to either
$\SU(2)\ltimes D_{2n+1}^\R$ or $\SU(2)\ltimes HU_{4n}^i$ for some
$i=0,1,2,3$ and some $n\ge 1$.\epr\label{SO3sg}
\end{prop}

% --------------------------------------------------------------------------------------------------------
%                                                                  section{application to Haag and Kaz}
%                                                                           subsection{reminder}
\section{Application to Haagerup and Kazhdan Properties}

\subsection{Reminder}\label{reminder}

Recall \cite[Chap. 1]{CCJJV} that a locally compact,
$\sigma$-compact group $G$ has the \textit{Haagerup Property} if
there exists a metrically proper, isometric action of $G$ on some
affine Hilbert space.

If $H$ is a subgroup of $G$, the pair $(G,H)$ has \textit{Kazhdan
Property (T)}, or that $H$ has Kazhdan's Property (T) relatively
to $G$, if every isometric action of $G$ on any affine Hilbert
space has a fixed point in restriction to $H$. In the case when
$H=G$, $G$ is said to have Property (T) (see \cite{HV} or
\cite{BHV}).

\medskip

As an immediate consequence of these definitions, if $(G,H)$ has
Property (T) and $H$ is not relatively compact in $G$, then $G$
does not have the Haagerup Property; this is a frequent
obstruction to Haagerup Property, although it is not the only one
(see Remark \ref{SansTrel_niHaag}).

The class of groups with the Haagerup Property generalizes the
class of amenable groups as a strong negation of Kazhdan's
Property (T). For other motivations of the Haagerup Property and
equivalent definitions, see \cite{CCJJV}.

In the following lemma, we summarize the hereditary properties of
the Haagerup and Kazhdan Properties that we will use in the
sequel.

\begin{lem}
The Haagerup Property for locally compact, $\sigma$-compact groups
is closed under taking (H1) closed subgroups, (H2) finite direct
products, (H3) direct limits \cite[Proposition 6.1.1]{CCJJV}, (H4)
extensions with amenable quotient \cite[Example 6.1.6]{CCJJV}, and
(H5) is inherited from lattices \cite[Proposition 6.1.5]{CCJJV}.

Relative Property (T) is inherited by dense images: if $(G,H)$ has
Property (T) and $f:G\to K$ is a continuous morphism, then
$(K,\overline{f(H)})$ has Property (T).\label{haagerup}
\end{lem}

\medskip

% --------------------------------------------------------------------------------------------------------
%                                                                         subsection{continuous families}

\subsection{Continuous families of Lie groups with Property (T)}~

\medskip

\noindent \textbf{Proof of Proposition \ref{cont_family_T1}.} We
must construct a continuous family of connected Lie groups with
Property (T) and with perfect and pairwise non-isomorphic Lie
algebras.

Consider $\mk{s}=\mk{sp}_{2n}(\R)$ ($n\ge 2$). Let $\mk{v}_i$,
$i=1,2,3,4$ be four nontrivial absolutely irreducible,
$\mk{s}$-modules which are pairwise non-isomorphic and all
preserve a symplectic form\footnote{There exist infinitely many
such modules, which can be obtained by taking large irreducible
components of the odd tensor powers of the standard
$2n$-dimensional $\mk{s}$-module.}. Then
$\mk{v}=\bigoplus_{i=1}^4\mk{v}_i$ is a full $\mk{s}$-module and
$\Aut_\mk{s}(\mk{v})=\prod_{i=1}^4\Aut_\mk{s}(\mk{v}_i)\simeq(\R^*)^4$.
In particular, $\Alt_\mk{s}(\mk{v})^*\simeq\R^4$ and
$\Aut_\mk{s}(\mk{v})$ acts diagonally on it. The action on the
2-Grassmannian, which is 4-dimensional, is trivial on the scalars,
so that its orbits are at most 3-dimensional. So there exists a
continuous family $(P_t)$ of 2-planes in $\Alt_\mk{s}(\mk{v})^*$
which are in pairwise distinct orbits for the action of
$\Aut_\mk{s}(\mk{v})$. By Theorem \ref{clasmin}, the Lie
$\mk{s}$-algebras $\mk{h}(\mk{v})/P_t$ are pairwise
non-isomorphic, and so the Lie algebras
$\mk{s}\ltimes\mk{h}(\mk{v})/P_t$ are pairwise non-isomorphic. The
Lie algebras $\g_t$ are perfect, and the corresponding Lie groups
$G_t$ have Property (T): this immediately follows from Wang's
classification \cite[Theorem 1.9]{Wang}.\epr

\begin{rem}
These examples have 2-nilpotent radical. This is, in a certain
sense, optimal, since there exist only countably many isomorphism
classes of Lie algebras over $\R$ with abelian radical, and only a
finite number for each dimension.
\end{rem}

\noindent \textbf{Proof of Proposition \ref{cont_family_T2}.} We
must construct a continuous family of locally isomorphic, pairwise
non-isomorphic connected Lie groups with Property (T). The proof
is actually similar to that of Proposition \ref{cont_family_T1}.
Use the same construction, but, instead of taking the quotient
$G_t$ by $P_t$, take the quotient $H_t$ by a lattice $\Gamma_t$ of
$P_t$. If we take the quotient of $H_t$ by its biggest compact
normal subgroup $P_t/\Gamma_t$, we obtain $G_t$. Accordingly, the
groups $H_t$ are pairwise non-isomorphic.\epr

\medskip

% --------------------------------------------------------------------------------------------------------
%                                                                           subsection{charact grps Haagerup}

\subsection{Characterization of groups with the Haagerup Property}

\begin{prop}
Let $\K$ be a local field of characteristic zero and $n\ge 1$.
Then the pairs $(\SL_2(\K)\ltimes V_n(\K),V_n(\K))$,
$(\PGL_2(\K)\ltimes V_n(\K),V_n(\K))$, $(\SL_2(\K)\ltimes
H_n(\K),H_n(\K))$, $(\widetilde{\SL_2(\R)}\ltimes
V_n(\R),V_n(\R))$, and $(\widetilde{\SL_2(\R)}\ltimes
H_n(\R),H_n(\R))$ have Property (T).\label{SL2relT}\end{prop}

\bpr The first (and the fourth) case is well-known; it follows,
for instance, from Furstenberg's theory \cite{FUR} of invariant
probabilities on projective spaces, which implies that $\SL_2(\K)$
does not preserve any probability on $V_n(\K)$ (more precisely, on
its Pontryagin dual) other than the Dirac measure at zero. See,
for instance, the proof of \cite[Chap. 2, Proposition 2]{HV}. The
second case is an immediate consequence of the first. For the
third (resp. fifth) case, we invoke \cite[Proposition
4.1.4]{CCJJV}, with $S=\SL_2(\K)$, $N=H_n(\K)$, even if the
hypotheses are slightly different (unless $\K=\R$ or $\C$): the
only modification is that, since here $[N,S]$ is not necessarily
connected, we must show that its image in the unitary group
$\textnormal{U}_n$ is connected so as to justify Lie's Theorem.
Otherwise, it would have a nontrivial finite quotient. This is a
contradiction, since $[N,S]$ is generated by divisible elements;
this is clear, since, as the group of $\K$-points of an unipotent
group, it has a well-defined logarithm.\epr

\begin{cor}
Let $G$ be either a connected Lie group, or
$G=\mathbf{G}(\mathbf{K})$, where $\mathbf{G}$ is a linear
algebraic group over the local field $\mathbf{K}$ of
characteristic zero. Suppose that the Lie algebra $\mathfrak{g}$
of $G$ contains a subalgebra $\mathfrak{h}$ isomorphic to either
$\mathfrak{sl}_2\ltimes \mathfrak{v}_n$ or $\mathfrak{sl}_2\ltimes
\mathfrak{h}_{2n-1}$ for some~$n\ge 2$. Then $G$ has a noncompact
closed subgroup with relative Property~(T). In particular, $G$
does not have Haagerup's property.\label{SncR=/1relT}
\end{cor}
\bpr Let us begin by the case of algebraic groups. By \cite[Chap.
II, Corollary 7.9]{Borel}, since $\mathfrak{h}$ is perfect, it is
the Lie algebra of a closed $\mathbf{K}$-subgroup $H$ of~$G$.
Since $H$ must be $\mathbf{K}$-isomorphic to either
$\textnormal{SL}_2\ltimes V_m$, $\textnormal{PGL}_2\ltimes
V_{2m-1}$, or $\textnormal{SL}_2\ltimes H_{2m-1}$ for some $m\ge
2$, Proposition \ref{SL2relT} implies that $G(\mathbf{K})$ has a
noncompact closed subgroup with relative Property~(T).

In the case of Lie groups, we obtain a Lie subgroup which is the
image of an immersion $i$ of
$\widetilde{\textnormal{SL}_2(\mathbf{R})}\ltimes N$, where $N$ is
either $V_n(\mathbf{R})$ or $H_{2n-1}(\mathbf{R})$, for some $n\ge
2$, into~$G$. By Proposition \ref{SL2relT}, $(G,\overline{i(N)})$
has Property~(T). We claim that $\overline{i(N)}$ is not compact.
Suppose the contrary. Then it is solvable and connected, hence it
is a torus. It is normal in the closure $H$ of~$i(G)$. Since the
automorphism group of a torus is totally disconnected, the action
by conjugation of $H$ on $\overline{i(N)}$ is trivial; that is,
$i(N)$ is central in~$H$. This is a contradiction.\epr

\medskip

\textbf{Proof of Theorem \ref{class_haag}.} As we already noticed
in the reminder, (i)$\Rightarrow$(ii) is immediate from the
definition. We are going to prove
(ii)$\Rightarrow$(iv)$\Rightarrow$(iii)$\Rightarrow$(i).

For the implication (iii)$\Rightarrow$(i), in the algebraic case,
$G$ is isomorphic, up to a finite kernel, to
$S_{nc}(\mathbf{K})\times \text{Mr}(\mathbf{K})$, where
$\text{Mr}$ denotes the amenable radical of~$\mathbf{G}$. The
group $\text{Mr}(\mathbf{K})$ is amenable, hence has Haagerup's
property. The group $S_{nc}(\mathbf{K})$ also has Haagerup's
property: if $\mathbf{K}$ is Archimedean, it maps, with finite
kernel, onto a product of groups isomorphic to
$\textnormal{PSO}_0(n,1)$ or $\textnormal{PSU}(n,1)$ ($n\ge 2$),
and these groups have Haagerup's property, by a result of Faraut
and Harzallah, see \cite[Chap. 2]{BHV}. If $\mathbf{K}$ is
non-Archimedean, then $S_{nc}(\mathbf{K})$ acts properly on a
product of trees (one for each simple factor) \cite{BT}, and this
also implies that it has Haagerup's property \cite[Chap. 2]{BHV}.

The same argument also works for connected Lie groups when the
semisimple part has finite centre; in particular, this is
fulfilled for linear Lie groups and their finite coverings. The
case when the semisimple part has infinite centre is considerably
more involved, see \cite[Chap. 4]{CCJJV}.

(ii)$\Rightarrow$(iv) Suppose that (iv) is not satisfied. If
$\mathfrak{g}$ contains a copy of $\mathfrak{sl}_2\ltimes
\mathfrak{v}_n$ or $\mathfrak{sl}_2\ltimes \mathfrak{h}_{2n-1}$
for some $n\ge 2$, then, by Corollary \ref{SncR=/1relT}, $G$ does
not satisfy (ii). If $\mathbf{K}=\mathbf{R}$, we consider $G$ as a
Lie group with finitely many components. By a standard argument,
since $\textnormal{Sp}(2,1)$ is simply connected with finite
centre (of order 2), an embedding of $\mathfrak{sp}(2,1)$ into
$\mathfrak{g}$ corresponds to a closed embedding of
$\textnormal{Sp}(2,1)$ or $\textnormal{PSp}(2,1)$ into~$G$. Since
$\textnormal{Sp}(2,1)$ has Property~(T) \cite[Chap. 3]{BHV}, this
contradicts~(ii).

(iv)$\Rightarrow$(iii) If $\mathfrak{g}$ is not M-decomposed,
then, by Proposition \ref{subalg_snc}, it contains a copy of
$\mathfrak{sl}_2\ltimes\mathfrak{v}_n$ or
$\mathfrak{sl}_2\ltimes\mathfrak{h}_{2n-1}$ for some~$n\ge 2$.

If $\mathfrak{g}$ has a simple factor $\mathfrak{s}$, then
$\mathfrak{s}$ embeds in $\mathfrak{g}$ through a Levi factor. If
$\mathfrak{s}$ has $\mathbf{K}$-rank $\ge 2$, then it contains a
subalgebra isomorphic to either $\mathfrak{sl_3}$ or
$\mathfrak{sp}_4$ \cite[Chap I, (1.6.2)]{Margulis}, and such a
subalgebra contains a subalgebra isomorphic to
$\mathfrak{sl}_2\ltimes\mathfrak{v}_2$ (resp.
$\mathfrak{sl}_2\ltimes\mathfrak{v}_3$) \cite[1.4 and 1.5]{BHV}.

Finally, if $\mathbf{K}=\mathbf{R}$ and $\mathfrak{s}$ is
isomorphic to either $\mathfrak{sp}(n,1)$ for some $n\ge 2$ or
$\mathfrak{f}_{4(-20)}$, then it contains a copy
of~$\mathfrak{sp}(2,1)$.\epr

\begin{rem}
Conversely, $\mk{sp}(n,1)$ does not contain any subalgebra
isomorphic to $\mk{sl}_2\ltimes\mk{v_n}$ or
$\mk{sl}_2\ltimes\mk{h_{2n-1}}$ for any $n\ge 2$; this can be
shown using results of \cite{CDSW} about weak amenability.
\end{rem}

\medskip

% --------------------------------------------------------------------------------------------------------
%                                                                           subsection{subgroups Lie groups}

\subsection{Subgroups of Lie groups}

Let us exhibit some subgroups in the groups above.

\begin{obs}
Let $G$ denote $\SL_2\ltimes V_n$, $\PGL_2\ltimes V_{2n-1}$, or
$\SL_2\ltimes H_{2n-1}$ for some $n\ge 2$, and $R$ its radical.
Then, for every field $K$ of characteristic zero, $G(K)$ contains
$G(\Z)$ as a subgroup. On the other hand, the pair $(G(\Z),R(\Z))$
has Property (T), this is because $G(\Z)$ is a lattice in $G(\R)$.
\end{obs}

\begin{obs}
Now, let $G$ denote $\SU(2)\ltimes D_{2n+1}^\R$, $\SO_3(\R)\ltimes
D_{2n+1}^\R$, or $\SU(2)\ltimes HU_{4n}^i$ for some $i=0,1,2,3$.
These groups are all defined over $\Q$: this is obvious at least
for all but $\SU(2)\ltimes HU_{4n}^i$ for $i=1,2$; for these two,
this is because the subspace $Z_i$ can be chosen rational in the
definition of $HU_{4n}^i$.

Let $R$ be the radical of $G$ and $S$ a Levi factor defined over
$\Q$. Let $F$ be a number field of degree three over $\Q$, not
totally real. Let $\mathcal{O}$ be its ring of integers. Then
$G(\mathcal{O})$ embeds diagonally as an irreducible lattice in
$G(\R)\times G(\C)$. Its projection $\Gamma$ in $G(\R)$ does not
have Haagerup's property, since otherwise $G(\C)$ would also have
Haagerup's property (by (H5) in Lemma \ref{haagerup}), and this is
excluded since it does not satisfy $[S_{nc},R]=1$, see Proposition
\ref{SncR=/1relT} (the anisotropic Levi factor becomes isotropic
after complexification). \label{obs2}\end{obs}

\begin{prop}
Let $G$ be a real Lie group, $R$ its radical, $S$ a semisimple
factor. Suppose that $[S,R]\neq 1$. Then $G$ has a countable
subgroup without Haagerup's property.\label{nh}
\end{prop}
\bpr First case: $[S_{nc},R]\neq 1$. Then, by Proposition
\ref{subliegrp_snc}, $G$ has a Lie subgroup $H$ isomorphic to a
quotient of $\widetilde{H}=\widetilde{\SL_2(\R)}\ltimes R(\R)$ by
a discrete central subgroup, where $R=V_n$ or $H_{2n-1}$, for some
$n\ge 2$. Denote by $\widetilde{H}(\Z)$ the inverse image of
$\SL_2(\Z)\ltimes R(\Z)$ in $\widetilde{H}$. By the observation
above, $(\widetilde{H}(\Z),R(\Z))$ has Property (T), so that its
image in $H$, which we denote by $H(\Z)$, satisfies
$(H(\Z),R_G(\Z))$ has Property (T), where $R_G(\Z)$ means the
image of $R(\Z)$ in $G$. Observe that $R_G(\Z)$ is infinite: if
$R=V_n$, this is $V_n(\Z)$; if $R=H_{2n-1}$, this is a quotient of
$H_{2n-1}(\Z)$ by some central subgroup. Accordingly, $H(\Z)$ does
not have Haagerup's property.

\medskip

Second case: $[S_c,R]\neq 1$. By Proposition \ref{SO3sg}, $G$ has
a Lie subgroup $H$ isomorphic to a central quotient of
$\SU(2)(\R)\ltimes R$, where $R=D_{2n+1}^\R$ or $HU_{4n}^i$, for
some $n\ge 1$ and $i=0,1,2,3$.

First suppose that the radical of $H$ is simply connected. Then,
by Observation \ref{obs2}, $H$ has a subgroup without the Haagerup
property.

Now, let us deal with the case when $H=\widetilde{H}/Z$, where $Z$
is a discrete central subgroup. Then $\widetilde{H}$ has a
subgroup $\Gamma$ as above which does not have Haagerup's
property. Let $W$ denote the centre of $\widetilde{H}$. The kernel
of the projection of $\Gamma$ to $H$ is given by $\Gamma\cap Z$.
We use the following trick: we apply an automorphism $\alpha$ of
$\widetilde{H}$ such that $\alpha(\Gamma)\cap Z$ is finite. It
follows that the image of $\alpha(\Gamma)$ in $H$ does not have
Haagerup's property.

This allows to suppose that $\Gamma\cap Z$ is finite, so that the
image of $\Gamma$ in $H$ does not have Haagerup's property. Let us
construct such an automorphism.

Observe that the representations of $\SU(2)$ can be extended to
the direct product $\R^*\times \SU(2)$ by making $\R^*$ act by
scalar multiplication. This action lifts to an action of
$\R^*\times \SU(2)$ on $HU_{4n}^i$, where the scalar $a$ acts on
the derived subgroup of $HU_{4n}^i$ by multiplication by $a^2$.

Now, working in the unit component of the centre $W$ of
$\widetilde{H}$, which we treat as a vector space, we can take $a$
so that $a^2\cdot(\Gamma\cap W)$ avoids $Z\cap W$ ($a$ clearly
exists, since $\Gamma$ and $Z$ are countable).\epr

\begin{defn}
Let $G$ be a locally compact group. We say that $G$ has Haagerup's
property if every $\sigma$-compact open subgroup of $G$ is.
\end{defn}

\begin{rem}
In view of (H3) of Lemma \ref{haagerup}, this is equivalent to:
every compactly generated, open subgroup of $G$ has Haagerup's
property, and also equivalent to the existence of a
$C_0$-representation with almost invariant vectors \cite[Chap.
1]{CCJJV}. In particular, $G$ having Haagerup's property and
$(G,H)$ having Property (T) still imply $H$ relatively compact.

All properties of the class of groups with Haagerup's property
claimed in Lemma \ref{haagerup} also clearly remain true for
general locally compact groups.
\end{rem}

If $G$ is a topological group, denote by $G_d$ the group $G$
endowed with the discrete topology.

\noindent \textbf{Proof of Theorem \ref{GdHaagerup}.} We remind
that we must prove, for a connected Lie group $G$, the equivalence
between

(i) $G$ is locally isomorphic to $\SO_3(\R)^a\times
\SL_2(\R)^b\times \SL_2(\C)^c\times R$, with $R$ solvable and
integers $a,b,c$, and

(ii) $G_d$ has Haagerup's property.

The implication (i)$\Rightarrow$(ii) is, essentially, a deep and
recent result of Guentner, Higson, and Weinberger \cite[Theorem
5.1]{GHW}, which implies that $(\PSL_2(\C))_d$ has Haagerup's
property. Let $G$ be as in (i), and $S$ its semisimple factor.
Then $G/S$ is solvable, so that, by (H4) of Lemma \ref{haagerup},
we can reduce to the case when $G=S$. Now, let $Z$ be the centre
of the semisimple group $G$, and embed $G_d$ in $(G/Z)_d\times G$,
where $G_d$ means $G$ endowed with the discrete topology. This is
a discrete embedding. Since $G$ has Haagerup's property, this
reduces the problem to the case when $G$ has trivial centre. So,
we are reduced to the cases of $\SO_3(\R)$, $\PSL_2(\R)$, and
$\PSL_2(\C)$. The two first groups are contained in the third, so
that the result follows from the Guentner-Higson-Weinberger
Theorem.

Conversely, suppose that $G$ does not satisfy (i).

If $[S,R]\neq 1$, then, by Proposition \ref{nh}, $G_d$ does not
have Haagerup's property. Otherwise, observe that the simple
factors allowed in (i) are exactly those of geometric rank one
(viewing $\SL_2(\C)$ as a complex Lie group). Hence, $S$ has a
factor $W$ which is not of geometric rank one. Then the result is
provided by Lemma \ref{higher_geom_rk} below.\epr

\begin{lem}
Let $S$ be a simple Lie group which is not locally isomorphic to
$\SO_3(\R)$, $\SL_2(\R)$ or $\SL_2(\C)$. Then $S_d$ does not have
Haagerup's property.\label{higher_geom_rk}
\end{lem}
\bpr Let $Z$ be the centre of $S$, so that $S/Z\simeq G(\R)$ for
some $\R$-algebraic group $G$. By assumption, $G(\C)$ has factors
of higher rank, hence does not have Haagerup's property. Let $F$
be a number field of degree three over $\Q$, not totally real. Let
$\mathcal{O}$ be its ring of integers. Then $G(\mathcal{O})$
embeds diagonally as an irreducible lattice in $G(\R)\times
G(\C)$, and is isomorphic to its projection in $G(\R)$. Let
$\Gamma$ be the inverse image in $S\times G(\C)$ of
$G(\mathcal{O})$. Then $\Gamma$ is a lattice in $S\times G(\C)$.
Hence, by \cite[Proposition 6.1.5]{CCJJV}, $\Gamma$ does not have
Haagerup's property. Note that the projection $\Gamma'$ of
$\Gamma$ into $S$ has finite kernel, contained in the centre of
$G(\C)$. So $\Gamma'$ neither has Haagerup's property, and is a
subgroup of $S$.\epr

\begin{rem}
Theorem \ref{GdHaagerup} is no longer true if we replace the
statement ``$G_d$ has Haagerup's property" by ``$G_d$ has no
infinite subgroup with relative Property (T)". Indeed, let
$G=K\ltimes V$, where $K$ is locally isomorphic to $\SO_3(\R)^n$
and $V$ is a vector space on which $K$ acts nontrivially. Suppose
that $(G_d,H)$ has Property (T) for some subgroup $H$. Then
$(G_d/V,H/(H\cap V))$ has Property (T). In view of the
Guentner-Higson-Weinberger Theorem (see the proof of Theorem
\ref{GdHaagerup}), $H/(H\cap V)$ is finite. On the other hand,
since $G$ has Haagerup's property, $H\cap V$ must be relatively
compact, and this implies that $H\cap V=1$. Thus, $H$ is finite.

Motivated by this example, it is easy to exhibit finitely
generated groups without the Haagerup Property and do not have
infinite subgroups with relative Property (T). For instance, let
$n\ge 3$, and $q$ be the quadratic form
$\sqrt{2}\,x_0^2+x_1^2+x_2^2+\dots+x_{n-1}^2$. Let
$G(R)=\SO(q)(R)\ltimes R^n$ and write, for any commutative
$\Q(\sqrt{2})$-algebra $R$, $H(R)=\SO(q)(R)$. Then
$\Gamma=G(\Z[\sqrt{2}])$ is such an example. The fact that
$\Gamma$ has no infinite subgroup $\Lambda$ with relative Property
(T) can be seen without making use of the
Guentner-Higson-Weinberger Theorem: first observe that
$H(\Z[\sqrt{2}])$ is a cocompact lattice in $\SO(n-1,1)$, hence
has Haagerup's property. So the projection of $\Lambda$ in
$H(\Z[\sqrt{2}])$ is finite. So, passing to a finite index
subgroup if necessary, we can suppose that $\Lambda$ is contained
in the subgroup $\Z[\sqrt{2}]^n$ of
$\Gamma=\SO(q)(\Z[\sqrt{2}])\ltimes \Z[\sqrt{2}]^n$. But then the
closure $L$ of $\Lambda$ in the subgroup $\R^n$ of the amenable
group $G(\R)=\SO(q)(\R)\ltimes \R^n$ is not compact, and
$(G(\R),L)$ has Property (T). This is a contradiction.

On the other hand, $\Gamma$ does not have Haagerup's property,
since it is a lattice in $G(\R)\ltimes G^\sigma(\R)$ (use (H5) of
Lemma \ref{haagerup}), where $\sigma$ is the nontrivial
automorphism of $\Q(\sqrt{2})$, and $G^\sigma(\R)\simeq
\SO(n-1,1)\ltimes\R^n$ does not have Haagerup's property, by
Corollary \ref{SncR=/1relT}. Note that $\Gamma$, as a cocompact
lattice in a connected Lie group, is finitely
presented.\label{SansTrel_niHaag}
\end{rem}

\textbf{Remerciements.} Je remercie Alain Valette pour m'avoir
lancé sur ces problèmes, pour ses nombreuses relectures des
versions préliminaires, et les corrections qui se sont ensuivies;
en particulier, je lui suis très reconnaissant de m'avoir
opportunément signalé la proposition 8.2 de \cite{CDSW}. Je
remercie également Pierre de la Harpe et le r\'ef\'er\'e anonyme
pour leurs judicieuses suggestions et corrections.

% --------------------------------------------------------------------------------------------------------
%                                                                             BIBLIOGRAPHIE

\bigskip

\footnotesize
\noindent Yves de Cornulier\\
\'Ecole Polytechnique Fédérale de Lausanne (EPFL)\\
Institut de Géométrie, Algèbre et Topologie (IGAT)\\
CH-1015 Lausanne, Switzerland\\
E-mail: \url{decornul@clipper.ens.fr}\\

% ----------------------------------------------------------------
\end{document}